\documentclass[11pt]{amsart}
\usepackage{amsopn}
\usepackage{amssymb, amscd}
\usepackage{multirow}
\usepackage{array}
\usepackage{graphicx, graphics, epsfig}
\usepackage{subcaption}
\usepackage{faktor}
\usepackage{color}
\usepackage{hyperref}
\usepackage{caption}
\usepackage{enumerate}
\usepackage{tikz-cd}
\usepackage{xcolor}

\graphicspath{ {imagenes/} }

\newcommand{\nc}{\newcommand}

\nc{\fg}{\mathfrak{f} } \nc{\vg}{\mathfrak{v} } \nc{\wg}{\mathfrak{w} }
\nc{\zg}{\mathfrak{z} } \nc{\ngo}{\mathfrak{n} } \nc{\kg}{\mathfrak{k} }
\nc{\mg}{\mathfrak{m} } \nc{\bg}{\mathfrak{b} } \nc{\ggo}{\mathfrak{g} } \nc{\eg}{\mathfrak{e} }
\nc{\ggob}{\overline{\mathfrak{g}} } \nc{\sog}{\mathfrak{so} }
\nc{\sug}{\mathfrak{su} } \nc{\spg}{\mathfrak{sp} } \nc{\slg}{\mathfrak{sl} }
\nc{\glg}{\mathfrak{gl} } \nc{\cg}{\mathfrak{c} } \nc{\rg}{\mathfrak{r} }
\nc{\hg}{\mathfrak{h} } \nc{\tg}{\mathfrak{t} } \nc{\ug}{\mathfrak{u} }
\nc{\dg}{\mathfrak{d} } \nc{\ag}{\mathfrak{a} } \nc{\pg}{\mathfrak{p} }
\nc{\sg}{\mathfrak{s} } \nc{\affg}{\mathfrak{aff} } \nc{\qg}{\mathfrak{q} } \nc{\lgo}{\mathfrak{l} }

\nc{\pca}{\mathcal{P}} \nc{\nca}{\mathcal{N}} \nc{\lca}{\mathcal{L}}
\nc{\oca}{\mathcal{O}} \nc{\mca}{\mathcal{M}} \nc{\tca}{\mathcal{T}}
\nc{\aca}{\mathcal{A}} \nc{\cca}{\mathcal{C}} \nc{\gca}{\mathcal{G}}
\nc{\sca}{\mathcal{S}} \nc{\hca}{\mathcal{H}} \nc{\bca}{\mathcal{B}}
\nc{\dca}{\mathcal{D}} \nc{\eca}{\mathcal{E}} \nc{\wca}{\mathcal{W}} \nc{\ica}{\mathcal{I}}

\nc{\vp}{\varphi} \nc{\ddt}{\tfrac{d}{dt}} \nc{\dsdt}{\tfrac{d^2}{dt^2}} \nc{\dds}{\frac{d}{ds}}
\nc{\dpar}{\frac{\partial}{\partial t}} \nc{\im}{\mathrm{i}}

\nc{\SO}{\mathrm{SO}} \nc{\Spe}{\mathrm{Sp}} \nc{\Sl}{\mathrm{SL}}
\nc{\SU}{\mathrm{SU}} \nc{\Or}{\mathrm{O}} \nc{\U}{\mathrm{U}} \nc{\Gl}{\mathrm{GL}}
\nc{\Se}{\mathrm{S}} \nc{\Cl}{\mathrm{Cl}} \nc{\Spin}{\mathrm{Spin}}
\nc{\Pin}{\mathrm{Pin}} \nc{\G}{\mathrm{GL}_n(\RR)} \nc{\g}{\mathfrak{gl}_n(\RR)}

\nc{\RR}{{\mathbb R}} \nc{\HH}{{\mathbb H}} \nc{\CC}{{\mathbb C}} \nc{\ZZ}{{\mathbb Z}}
\nc{\FF}{{\mathbb F}} \nc{\NN}{{\mathbb N}} \nc{\QQ}{{\mathbb Q}} \nc{\PP}{{\mathbb P}} \nc{\OO}{{\Bbb O}}

\nc{\vs}{\vspace{.2cm}} \nc{\vsp}{\vspace{1cm}} \nc{\ip}{\langle\cdot,\cdot\rangle}
\nc{\ipp}{(\cdot,\cdot)} \nc{\la}{\langle} \nc{\ra}{\rangle} \nc{\unm}{\tfrac{1}{2}}
\nc{\unc}{\tfrac{1}{4}} \nc{\und}{\frac{1}{16}} \nc{\no}{\vs\noindent}
\nc{\lam}{\Lambda^2(\RR^n)^*\otimes\RR^n} \nc{\tangz}{{\rm T}^{\rm Zar}}
\nc{\nor}{{\sf n}}  \nc{\mum}{/\!\!/} \nc{\kir}{/\!\!/\!\!/}
\nc{\Ri}{\tfrac{4\Ric_{\mu}}{||\mu||^2}} \nc{\ds}{\displaystyle}
\nc{\ben}{\begin{enumerate}} \nc{\een}{\end{enumerate}} \nc{\f}{\frac}
\nc{\lb}{[\cdot,\cdot]} \nc{\isn}{\tfrac{1}{||v||^2}}
\nc{\gkp}{(\ggo=\kg\oplus\pg,\ip)} \nc{\ukh}{(\ug=\kg\oplus\hg,\ip)}
\nc{\tgkp}{(\tilde{\ggo}=\kg\oplus\pg,\ip)}
\nc{\wt}{\widetilde}
\nc{\iop}{\mathtt{i}} \nc{\jop}{\mathtt{j}}
\nc{\Hk}{H_{\kil}} \nc{\gk}{g_{\kil}}

\nc{\Hess}{\operatorname{Hess}} \nc{\ad}{\operatorname{ad}}
\nc{\Ad}{\operatorname{Ad}} \nc{\rank}{\operatorname{rk}}
\nc{\Irr}{\operatorname{Irr}} \nc{\End}{\operatorname{End}}
\nc{\Aut}{\operatorname{Aut}} \nc{\Inn}{\operatorname{Inn}}
\nc{\Der}{\operatorname{Der}} \nc{\Ker}{\operatorname{Ker}}
\nc{\Iso}{\operatorname{Iso}} \nc{\Diff}{\operatorname{Diff}}
\nc{\Lie}{\operatorname{L}} \nc{\tr}{\operatorname{tr}} \nc{\dif}{\operatorname{d}}
\nc{\sen}{\operatorname{sen}} \nc{\modu}{\operatorname{mod}}
\nc{\CRic}{\operatorname{PP}} \nc{\Cric}{\operatorname{P}} \nc{\Ricci}{\operatorname{Ric}}
\nc{\sym}{\operatorname{sym}} \nc{\herm}{\operatorname{herm}} \nc{\symac}{\operatorname{sym^{ac}}}
\nc{\symc}{\operatorname{sym^{c}}} \nc{\scalar}{\operatorname{Sc}}
\nc{\grad}{\operatorname{grad}} \nc{\ricci}{\operatorname{Rc}} \nc{\kil}{\operatorname{B}} \nc{\cas}{\operatorname{C}} \nc{\lic}{\operatorname{L}}
\nc{\Nor}{\operatorname{Norm}}  \nc{\ricc}{\operatorname{Rc^{c}}}
\nc{\Ricc}{\operatorname{Ric^{c}}} \nc{\ricac}{\operatorname{Rc^{ac}}}
\nc{\Ricac}{\operatorname{Ric^{ac}}} \nc{\Riem}{\operatorname{Rm}} \nc{\Sec}{\operatorname{Sec}}
\nc{\riccig}{\operatorname{ric^{\gamma}}} \nc{\mm}{\operatorname{m}} \nc{\Mm}{\operatorname{M}}
\nc{\Le}{\operatorname{L}} \nc{\tang}{\operatorname{T}}
\nc{\level}{\operatorname{level}} \nc{\rad}{\operatorname{r}}
\nc{\abel}{\operatorname{ab}} \nc{\CH}{\operatorname{CH}} \nc{\Cone}{{\mathcal C}} \nc{\CCone}{\operatorname{CC}} \nc{\CP}{{\mathcal P}}
\nc{\mcc}{\operatorname{mcc}} \nc{\Adj}{\operatorname{Adj}}
\nc{\Order}{\operatorname{O}}  \nc{\inj}{\operatorname{inj}} \nc{\proy}{\operatorname{pr}}
\nc{\vol}{\operatorname{vol}} \nc{\Diag}{\operatorname{Dg}} \nc{\Diagg}{\operatorname{Diag}}
\nc{\Spec}{\operatorname{Spec}} \nc{\Ima}{\operatorname{Im}} \nc{\Rea}{\operatorname{Re}}
\nc{\spann}{\operatorname{span}} \nc{\Aff}{\operatorname{Aff}} \nc{\E}{\operatorname{E}} \nc{\id}{\operatorname{id}} \nc{\dete}{\operatorname{det}} \nc{\Crit}{\operatorname{Crit}} \nc{\val}{\operatorname{val}}

\theoremstyle{plain}
\newtheorem{theorem}{Theorem}[section]
\newtheorem{proposition}[theorem]{Proposition}
\newtheorem{corollary}[theorem]{Corollary}

\theoremstyle{definition}

\theoremstyle{remark}
\newtheorem{remark}[theorem]{Remark}

\title{Stability of non-diagonal Einstein metrics on homogeneous spaces $H\times H/\Delta K$}

\author{Valeria Guti\'errez}
\address{FAMAF, Universidad Nacional de C\'ordoba and CIEM, CONICET (Argentina)}
\email{valeria.gutierrez@unc.edu.ar}

\begin{document}

\begin{abstract}
We consider the homogeneous space $M=H\times H/\Delta K$, where $H/K$ is an irreducible symmetric space and $\Delta K$ denotes diagonal embedding. Recently, Lauret and Will provided a complete classification of $H\times H$-invariant Einstein metrics on M. They obtained that there is always at least one non-diagonal Einstein metric on $M$, and in some cases, diagonal Einstein metrics also exist. We give a formula for the scalar curvature of a subset of $H\times H$-invariant metrics and study the stability of non-diagonal Einstein metrics on $M$ with respect to the Hilbert action, obtaining that these metrics are unstable with different coindexes for all homogeneous spaces $M$.
\end{abstract}

\maketitle
\tableofcontents

\section{Introduction}
The known results on the existence and classification of Einstein metrics on compact homogeneous spaces $G/K$ vary significantly depending on whether $G$ is simple or not. Indeed, when $G$ is not simple, the isotropy representation of $G/K$ is almost never multiplicity-free and much more complicated. Recently, Lauret and Will on \cite{HHK} provided a complete classification of $H\times H$-invariant Einstein metric on any compact homogeneous space of the form
$$
M=H\times H/\Delta K,
$$
where $\Delta K:=\{(k,k)\in H\times H :k\in K\}$ and $H/K$ is an irreducible symmetric space such that $\kil_\kg=a\kil_\hg|_\kg$, where $\kil_\kg$ and $\kil_\hg$ are the Killing forms of Lie algebras $\kg$ and $\hg$ of $K$ and $H$, respectively. They also found an Einstein metric on $M$ for every irreducible symmetric space $H/K$.

If $G:=H\times H$ and $\hg=\kg\oplus\qg$ is the $\kil_{\hg}$-orthogonal reductive decomposition for $H/K$, then the $\kil_{\ggo}$-orthogonal reductive decomposition for $M=G/\Delta K$ is given by $\ggo=\hg \oplus \hg = \Delta\kg\oplus\pg$, where $\Delta\kg=\{(Z,Z):Z\in\kg\}$ and
$$
\pg=\pg_1\oplus\pg_2\oplus\pg_3,
\quad
\pg_1=(\qg,0), \quad \pg_2=(0,\qg), \quad \pg_3=\{(Z,-Z):Z\in\kg\}.
$$

If $\mca^G$ is the space of all $G$-invariant metrics on $M$, then we consider the subspace $\mca^{4}\subset \mca^G$ of metrics $g=(x_1,x_2,x_3,x_4)$ on $M$  defined as follows:
$$
g|_{\pg_i\times\pg_i}=x_i(-\kil_\ggo)|_{\pg_i\times\pg_i}, \; i=1,2,3, \quad g|_{\pg_1\times\pg_2}=x_4(-\kil_\ggo)|_{\pg_1\times\pg_2}, \quad g|_{\pg_1\times\pg_3}=g|_{\pg_2\times\pg_3}=0,
$$
where $x_1,x_2,x_3>0$ and $x_1x_2>x_4^2$. These metrics are called {\it diagonal} when $x_4=0$ and $\mca^{4} = \mca^G$ if and only if $K$ is either simple or one-dimensional.

In this paper, we provide a formula for the scalar curvature $\scalar(g)=\scalar(x_1,x_2,x_3,x_4)$ for any $g \in \mca^4$ (see Proposition \ref{scalar-fun}), which is given by
$$\scalar(g)=\tfrac{d}{2x_3}-\tfrac{n \left(x_2^2 x_3^2 + x_1^2 (4 x_2^2 - 8 x_2 x_3 + x_3^2) + 8 x_2 x_3 x_4^2 +
8 x_1 x_3 (x_4^2-x_2^2) - 2 x_4^2 (x_3^2 + 2 x_4^2)\right)}{16x_3(x_4^2-x_1 x_2)^2},$$
where $d:=\dim{K}$ and $n:=\dim{H/K}$.

In order to study the stability types of the non-diagonal Einstein metrics as critical points of the scalar curvature functional, we consider
$$
\scalar:\mca^{4}_1 \longrightarrow \RR,
$$
where $\mca^{4}_1$ is the manifold of all unit volume $G$-invariant metrics on $\mca^{4}$. An Einstein metric $g\in\mca_1^{4}$ is called $G$-{\it unstable} if $\scalar''_g(T,T)>0$ for some $T\in{T_g\mca_1^{4}}$, and is called \textit{non-degenerate} if $\scalar''_g|_{T_g\mca_1^{4}}$ is non-degenerate.

In \cite[Theorem 7.3]{HHK} (see Theorem \ref{non-diag}), Lauret and Will proved that if $\kil_\kg=a\kil_\hg|_\kg$ and $a>\tfrac{1}{2}$, then there exist two non-homothetic, non-diagonal $H\times H$-invariant Einstein metrics on the homogeneous space $M$. These metrics are given by $g_3=(1,1,1,y_+)$, where $y_+=\frac{1}{2}\sqrt{\frac{2a-1}{2-a}}$, and $g_5=\left(\unm,\frac{3}{2},1,\unm\right)$. Actually, $g_5$ is also Einstein on $M$ for any symmetric space $H/K$. When $a<\tfrac{1}{2}$, they also found diagonal ($x_4=0$) Einsten metrics on $M$, the stability of these metrics was studied in \cite[Section 6]{HHK},  where it was demonstrated that these metrics are both $G$-unstable with different coindexes.

Our main result is as follows:
\begin{theorem}
For any irreducible symmetric space $H/K$, the non-diagonal Einstein metrics $g_3$ and $g_5$ on $M=H\times H/\Delta K$ are non-degenerate, $G$-unstable, with coindexes equal to $2$ and $1$, respectively (see Figure \ref{figu}).
\end{theorem}

The Einstein metric $g_5$ is most likely the one predicted by the Graph Theorem (see \cite[Theorem C]{BhmWngZll}). The class of spaces involved in the above theorem consists of $3$ infinite families and $2$ sporadic examples such that $\kil_\kg=a\kil_\hg|_\kg$ for $a>\tfrac{1}{2}$ (see \cite[Appendix A. Table 3]{HHK}), as well as $5$ more families and $7$ isolated examples where only $g_5$ exists (see Table \ref{irr-spaces}).

\vs \noindent \textit {Acknowledgements.} I would like to thank to my Ph.D. advisor Dr. Jorge Lauret for his
continued guidance during the preparation of this paper.

\section{Preliminaries}
\subsection{Homogeneous spaces $H\times H / \Delta K$}
The known results on the existence and classi\-fication of Einstein metrics on compact homogeneous spaces $G/K$ differs substantially between the cases of $G$ simple and non-simple. Recently, Lauret and Will in \cite{HHK} found new examples of invariant Einstein metrics on spaces of the form $H\times H/\Delta K$ with $H$ simple and $\Delta K$ denoting diagonal embedding.

Given a homogeneous space $H/K$ of dimension $n$, where $H$ is a compact simple Lie group and $K$ a proper closed subgroup of it with $d:=\dim{K}>0$, we consider the $(2n+d)$-dimensional compact homogeneous space
$$
M^{2n+d}=H\times H/\Delta K,
\qquad \text{ where } \quad \Delta K:=\{(k,k)\in H\times H :k\in K\}.
$$
Let $\hg=\kg\oplus\qg$ be the $\kil_{\hg}$-orthogonal reductive decomposition for $H/K$, where $\hg$, $\kg$ denote the corresponding Lie algebras and $\kil_{\hg}$ denotes the Killing form of $\hg$. If we set $G:=H\times H$ and $\ggo=\hg\oplus\hg$ its Lie algebra, then we have the $\kil_{\ggo}$-orthogonal reductive decomposition for $M$ given by
$$\ggo=\Delta\kg\oplus\pg \qquad \text{ where } \quad \Delta\kg=\{(Z,Z):Z\in\kg\},$$
and
\begin{equation}\label{pp}
  \pg=\pg_1\oplus\pg_2\oplus\pg_3,
\quad
\pg_1=(\qg,0), \quad \pg_2=(0,\qg), \quad \pg_3=\{(Z,-Z):Z\in\kg\}.
\end{equation}
The decomposition $\pg=\pg_1\oplus\pg_2\oplus\pg_3$ is $\Ad(K)$-invariant and $\gk$-orthogonal, where $\gk$ is the standard metric on $M=G/\Delta K$. As $\Ad(K)$-representations, $\pg_1\simeq\pg_2\simeq\qg$, the isotropy representation of $H/K$ (so $G/\Delta K$ is never multiplicity-free), and $\pg_3\simeq\kg$, the adjoint representation of $K$. We also consider the $\gk$-orthogonal, $\Ad(K)$-invariant decomposition
$$
\pg_3=\pg_3^0\oplus\dots\oplus\pg_3^t, \qquad \pg_3^l=\{(Z,-Z):Z\in\kg_l\},
$$
where
\begin{equation}\label{decs}
\kg=\kg_0\oplus\kg_1\oplus\dots\oplus\kg_t,
\end{equation}
$\kg_0$ is the center of $\kg$ and $\kg_1,\dots\kg_t$ are the simple ideals of $\kg$.

The equivalence of $\pg_1$ and $\pg_2$ gives rise to $G$-invariant metrics on the homogeneous space $M$ of the form
$$
g|_{\pg_i\times\pg_i}=x_i(-\kil_\ggo)|_{\pg_i\times\pg_i}, \; i=1,2,3, \quad g|_{\pg_1\times\pg_2}=x_4(-\kil_\ggo)|_{\pg_1\times\pg_2}, \quad g|_{\pg_1\times\pg_3}=g|_{\pg_2\times\pg_3}=0,
$$
where $x_1,x_2,x_3>0$ and $x_1x_2>x_4^2$. These metrics are denoted by
\begin{equation}\label{g-nondiag}
g=(x_1,x_2,x_3,x_4)
\end{equation}
and called {\it diagonal} when $x_4=0$. We denote by $\mca^{4}$ the space of all $G$-invariant metrics on $M$ of this form and note that if the isotropy representation of $H/K$ is irreducible and of real type and $K$ is either simple or one-dimensional, then $\mca^{4}=\mca^G$ the space of all $G$-invariant metrics on $M$.

\subsection{Non-diagonal Einstein Metrics on $H\times H / \Delta K $}\label{metric}
In \cite{HHK}, conditions for the existence of diagonal Einstein metrics on spaces of the form $M=H\times H / \Delta K $ were provided. Remarkably, most non-existence cases arose when $H/K$ is an irreducible symmetric space. However, it is also proved in \cite[Section 7]{HHK} that in these cases, there always exists a non-diagonal Einstein metric on $M$.

Consider a $G$-invariant metric $g=\left(x_1,x_2,x_3,x_4\right)$ as in \eqref{g-nondiag}. Note that for all $X_i,Y_i\in\qg$, $Z,W\in\kg$,
\begin{align}
g((X_1,X_2),(Y_1,Y_2)):=& x_1(-\kil_\hg)(X_1,Y_1) +x_2(-\kil_\hg)(X_2,Y_2) \notag \\
&+x_4\left((-\kil_\hg)(X_1,Y_2)+(-\kil_\hg)(X_2,Y_1)\right),  \label{defg4} \\
g((Z,-Z),(W,-W)):=& x_3(-\kil_\ggo)((Z,-Z),(W,-W)) = 2x_3(-\kil_\hg)(Z,W), \notag
\end{align}
and the matrix of $g$ with respect to $\gk$, the standard metric on $M$, is therefore given by
\begin{equation}\label{mat}
\left[\begin{matrix}
x_1I_{n}&x_4I_{n}&0\\
x_4I_{n}&x_2I_{n}&0\\
0&0&x_3I_{d}
\end{matrix}\right],
\end{equation}
where $d:=\dim{K}$ and $n:=\dim{H/K}$.

If we define, $$
\begin{array}{c}
\qg_1:=\pg_1=\left\{(X,0):X\in\qg\right\}, \qquad \qg_2:=\left\{(-x_4X,x_1X):X\in\qg\right\}, \\ \\ \qg_3:=\pg_3=\left\{(Z,-Z):Z\in\kg\right\}.
\end{array}
$$
then,
\begin{equation}\label{ort-dec}
\pg=\qg_1\oplus\qg_2\oplus\qg_3,
\end{equation}
is a $g$-orthogonal $\Ad(K)$-invariant decomposition of $\pg$.

An explicit description of the Ricci tensor for these non-diagonal metrics when $H/K$ is an irreducible symmetric space was given in \cite[Section 7]{HHK}.

\begin{proposition}{\upshape\cite[Proposition 7.1]{HHK}}\label{ric-nondiag}
If $H/K$ is an irreducible  symmetric space, then the Ricci tensor of the metric $g=(x_1,x_2,x_3,x_4)$ on $M^{2n+d}=H\times H/\Delta K$ is given for any
$$
\overline{X}=(X,0)\in\pg_1, \quad X\in\qg, \quad \overline{Y}=(0,Y)\in\pg_2, \quad Y\in\qg, \quad \overline{Z}=(Z,-Z)\in\pg_3^l, \quad Z\in\kg_l,
$$
$l=0,\dots,t$, as follows: $\quad\ricci(g)(\overline{X},\overline{Z}) = 0$, $\quad \ricci(g)(\overline{Y},\overline{Z}) = 0$,
\begin{align*}
\ricci(g)(\overline{X},\overline{X}) =& \left(\frac{x_3}{8x_1} + \frac{x_3x_4^2}{8x_1(x_1x_2-x_4^2)}
-\frac{x_1x_4^2}{2(x_1x_2-x_4^2)x_3} - \frac{1}{2}\right) \kil_\hg(X,X),\\
\ricci(g)(\overline{Y},\overline{Y}) =& \left(\frac{x_1x_3}{8(x_1x_2-x_4^2)}
+\frac{x_1x_2-x_4^2}{8x_1x_3}
-\frac{(x_1x_2+x_4^2)^2}{8x_1(x_1x_2-x_4^2)x_3} -  \frac{1}{2}\right) \kil_\hg(Y,Y),\\
\ricci(g)(\overline{X},\overline{Y}) =& \left(\frac{x_3x_4}{8(x_1x_2-x_4^2)}
-\frac{x_4}{4x_3}
- \frac{x_4(x_1x_2+x_4^2)}{4(x_1x_2-x_4^2)x_3}\right)  \kil_\hg(X,Y),\\
\ricci(g)(\overline{Z},\overline{Z}) =& \left(a_l(R-1)-R\right)
\kil_\hg(Z,Z), \qquad \kil_{\kg_l}=a_l\kil_\hg|_{\kg_l},
\end{align*}
where
\begin{equation}\label{defR}
R:=
-\frac{2x_4^2}{x_1x_2-x_4^2}
+ \frac{x_3^2}{4x_1^2}
+\frac{x_3^2(x_1^2-x_4^2)^2}{4x_1^2(x_1x_2-x_4^2)^2}
+\frac{x_3^2x_4^2}{2x_1^2(x_1x_2-x_4^2)}.
\end{equation}
\end{proposition}

The non-diagonal context provides many Einstein metrics, including one which exists on any $M=H\times H/\Delta K$, where $H/K$ is an irreducible symmetric space.

\begin{theorem}{\upshape\cite[Theorem 7.3]{HHK}}\label{non-diag} 
Let $H/K$ be an irreducible symmetric space.
\begin{enumerate}[{\rm (i)}]
\item If $\kil_\kg=a\kil_\hg|_\kg$, $a>0$, then, up to scaling, the $H\times H$-invariant Einstein metrics on the homogeneous space $M=H\times H/\Delta K$ are given by
\begin{enumerate}[{\small $\bullet$}]
\item $g_1=(x_+,x_+,1,0)$ and $g_2=(x_-,x_-,1,0)$ if $a<\unm$, where
$$
x_\pm=\frac{1\pm \sqrt{1 - a(3-2a)}}{2a}.
$$
\item $g_3=(1,1,1,y_+)$ and $g_4=(1,1,1,y_-)$ if $a>\unm$, where
$$
y_\pm=\pm\frac{1}{2}\sqrt{\frac{2a-1}{2-a}}.
$$
\item $g_5=\left(\unm,\frac{3}{2},1,\unm\right)$ and $g_6=\left(\frac{3}{2},\unm,1,\unm\right)$.
\end{enumerate}

\item The normalized scalar curvatures $\scalar_N(g_i):=\scalar(g_i)(\det_{\gk}{g_i})^{\frac{1}{\dim{M}}}$ are given as follows:
$$
\scalar_N(g_1)=\frac{(2n+d)(4x_+-1)}{8(x_+)^{2\alpha}}, \qquad
\scalar_N(g_2)=\frac{(2n+d)(4x_--1)}{8(x_-)^{2\alpha}},
$$
$$
\scalar_N(g_3)=\scalar_N(g_4)=\frac{3(2n+d)}{8} \left(\frac{4(2-a)}{3(3-2a)}\right)^\alpha, \qquad
\scalar_N(g_5)=\scalar_N(g_6)=(2n+d)2^{\alpha-2},
$$
where $\alpha=\tfrac{n+d}{2n+d}$.

\item The metrics $g_5$ and $g_6$ are also Einstein on $M=H\times H/\Delta K$ for any $H/K$.
\end{enumerate}
\end{theorem}

Note that the metrics $g_5$ and $g_6$ are isometric via the automorphism of $G=H\times H$ that exchanges the two factors of $H$, and considering $\theta$ as the idempotent automorphism of $H$ such that $d\theta|_\kg=I$ and $d\theta|_\qg=-I$, we obtain that $(id,\theta)\in\Aut(G/\Delta K)$ defines an isometry between $g_3$ and $g_4$ (\cite[Remark 7.6]{HHK}).

In \cite[Appendix A. Table 3]{HHK}, the irreducible symmetric spaces $H/K$ such that the condition $\kil_\kg=a\kil_\hg|_\kg$ holds are listed. When $a>\tfrac{1}{2}$, Theorem \ref{non-diag} states that there are two non-diagonal Einstein metrics on $M$, up to isometry and scaling. Conversely, if $a<\tfrac{1}{2}$ or if $H/K$ is among the other irreducible symmetric spaces given in Table \ref{irr-spaces} (see also \cite[Table 7.102]{Bss}), there exists a unique non-diagonal Einstein metric, up to isometry, given by $g_5=(\tfrac{1}{2},\tfrac{3}{2},1,\tfrac{1}{2}).$

\begin{table}[ht]
{\small
$$
\begin{array}{c|c|c|c|c}
H/K & m,p \text{ or } q & \dim{H} & d=\dim{K} & n=\dim{H/K}  
\\[2mm] \hline \hline \rule{0pt}{14pt}
\text{\footnotesize $\SU(p+q)/\SU(p)\times \SU(q)\times S^1$} & 1\leq p \leq q & (p+q)^2-1 & p^2+q^2-1 & 2pq
\\[2mm]  \hline \rule{0pt}{14pt}
\SO(2m)/\U(m) & m \geq 2 & 2m^2-m & m^2 & m(m-1)
\\[2mm]  \hline \rule{0pt}{14pt}
\text{\footnotesize $\SO(p+q)/\SO(p) \times \SO(q)$} & \substack{1 \ \leq \ p \ < \ q, \\[1.3mm] p+q \ \geq \ 7} & \tfrac{p^2+q^2+2pq-p-q}{2} & \tfrac{p^2-p+q^2-q}{2} & pq
\\[2mm]  \hline \rule{0pt}{14pt}
\text{\footnotesize $\Spe(m)/\SU(m)\times S^1$} & m \geq 2 & m(2m+1)& m^2 & m(m+1)
\\[2mm]  \hline \rule{0pt}{14pt}
\text{\footnotesize $\Spe(p+q)/\Spe(p)\times\Spe(q)$}  & 1\leq p <q & \text{\tiny $(p+q)(2p+2q+1)$} & \text{\tiny $p(2p+1)+q(2q+1)$} & 4pq
\\[2mm]  \hline \rule{0pt}{14pt}
G_2/\SO(4) &  & 14 & 6 & 8 
\\[2mm]  \hline \rule{0pt}{14pt}
F_4/\Spe(3)\times \SU(2) &  & 52 & 24 & 28 
\\[2mm]  \hline \rule{0pt}{14pt}
E_6/\SU(6) \times \SU(2) &  & 78 & 38 & 40 
\\[2mm]  \hline \rule{0pt}{14pt}
E_6/\SO(10) \times S^1 &  & 78 & 46 & 32 
\\[2mm]  \hline \rule{0pt}{14pt}
E_7/\SO(12)\times \SU(2) &  & 133 & 69 & 64 
\\[2mm] \hline  \rule{0pt}{14pt}
E_7/E_6\times S^1 &  & 133 & 79 & 54 
\\[2mm] \hline  \rule{0pt}{14pt}
E_8/E_7\times \SU(2) &  & 248 & 139 & 109 
\\[2mm] \hline\hline
\end{array}
$$}
\caption{Irreducible symmetric spaces $H/K$ such that $\kil_\kg\neq a\kil_\hg|_\kg$, for all $a\in \RR$. Those such that $\kil_\kg=a\kil_\hg|_\kg$, for some $a>0$ are listed in \cite[Appendix A. Table 3]{HHK}.}
\label{irr-spaces}.
\end{table}

\section{Stability}

In this section, we study the behavior of the scalar curvature function restricted to $\mca_1^{4}$,
$$\scalar: \mca^{4}_1 \rightarrow \RR,$$
where $\mca_1^{4}\subset \mca^{4}$ is the space of all $G$-invariant unit volume metrics defined as in (\ref{g-nondiag}) on $M=G/\Delta K$.

According to Hilbert \cite{Hlb}, the critical points of this functional are precisely the Einstein metrics of this form on $M$. Based on this, if $\scalar''_g(T,T)>0$ for an Einstein metric $g\in\mca_1^{4}$ and some $T\in T_g\mca_1^{4}$, then $g$ is called  $G$-{\it unstable}, and if $\scalar''_g|_{T_g\mca_1^{4}}$ is non-degenerate, then $g$ is called \textit{non-degenerate}.

In \cite[Section 6]{HHK}, it was studied the stability of diagonal Einstein metrics on $M=H\times H/\Delta K$ (i.e. $x_4=0$), and it was proved that these metrics are both $G$-unstable with different coindexes.

\begin{proposition}\label{scalar-fun}
If $H/K$ is an irreducible  symmetric space, the scalar curvature of any metric of the form $g=(x_1,x_2,x_3,x_4)$ on $M=H\times H / \Delta K$ is given by:
$$\scalar(g)=\tfrac{d}{2x_3}-\tfrac{n \left(x_2^2 x_3^2 + x_1^2 (4 x_2^2 - 8 x_2 x_3 + x_3^2) + 8 x_2 x_3 x_4^2 +
   8 x_1 x_3 (x_4^2-x_2^2) - 2 x_4^2 (x_3^2 + 2 x_4^2)\right)}{16x_3(x_4^2-x_1 x_2)^2}.$$
\end{proposition}
\begin{remark}
  If $K$ is simple or one-dimensional, this is the scalar curvature of all $G$-invariant metrics on $M$.
\end{remark}
\begin{remark}
  This formula is consistent with the one given in  \cite[Lemma 4.4]{HHK} for normal metrics $g_b=(z_1,z_2,\tfrac{2z_1z_2}{z_1+z_2},0)$ when $H/K$ is an irreducible symmetric space (i.e. $S=d-\tfrac{n}{2}$).
\end{remark}
\begin{proof}

Since $\scalar(g)=\tr_g\ricci(g)$, we consider the $g$-orthogonal decomposition $$\pg=\qg_1\oplus\qg_2\oplus\qg_3,$$ (see \eqref{ort-dec}), and the $g$-orthonormal bases of $\qg_1$, $\qg_2$ and $\qg_3$, respectively, given by:
$$
\left\{Y^1_\alpha:=\tfrac{1}{\sqrt{x_1}}(e_\alpha,0)\right\}, \quad
\left\{Y^2_\alpha:=\tfrac{1}{\sqrt{x_1(x_1x_2-x_4^2)}}(-x_4e_\alpha,x_1e_\alpha)\right\}, \quad
\left\{ Y^3_\alpha:=\tfrac{1}{\sqrt{2x_3}}(Z_\alpha,-Z_\alpha)\right\},
$$
where $\{ Z_\alpha\}$ is a $(-\kil_\hg)$-orthonormal basis of $\kg$, adapted to the decomposition \eqref{decs}, and $\{ e_\alpha\}$ is a $(-\kil_\hg)$-orthonormal basis of $\qg$ (see \cite[Section 7]{HHK}).

It is easy to see that,
\begin{align*}
\ricci(g)(Y^1_\alpha,Y^1_\alpha) =& -\tfrac{1}{x_1}\ricci(g)((e_\alpha,0),(e_\alpha,0))\\
\ricci(g)(Y^2_\alpha,Y^2_\alpha)=&\tfrac{1}{x_1(x_1x_2-x_4^2)}\left(-x_4^2\ricci(g)((e_\alpha,0),(e_\alpha,0))\right.\\
& \left. \qquad \qquad +2x_4x_1\ricci(g)((e_\alpha,0),(0,e_\alpha))-x_1^2\ricci(g)\left((0,e_\alpha),(0,e_\alpha)\right)\right)\\
\ricci(g)(Y^3_\alpha,Y^3_\alpha) =& \tfrac{1}{2x_3}\left(a_l(1-R)+R\right),\\
\end{align*}
where $R$ is as in \eqref{defR}. Then the proposition follows from Proposition \ref{ric-nondiag} and the fact that
\begin{equation}\label{trcas}
\tr{\cas_{\chi}}=\sum_{l=0}^t (1-a_l)d_l=\tfrac{n}{2}
\end{equation}
for a $n$-dimensional symmetric space $H/K$, where $\cas_{\chi}$ is the Casimir operator of the isotropy representation of $H/K$ with respect to $-\kil_{\hg}|_\kg$.
\end{proof}

The formula for the normalized scalar curvature of any $G$-invariant metric of the form $g=(x_1,x_2,1,x_4)$ is provided in the following proposition. It is worth noting that for Einstein metrics of this type, the formula was obtained in \cite[Theorem 7.3]{HHK} (see Theorem \ref{non-diag}).

\begin{corollary}\label{snor}
 The normalized scalar curvature $\scalar_N(g):=\scalar(g)(\det_{\gk}{g})^{\frac{1}{\dim{M}}}$ for any metric of the form $g=(x_1,x_2,1,x_4)$ on $M = H \times H/\Delta K$, where $H/K$ is an irreducible symmetric space is given by:
$$\scalar_N(g)=\tfrac{(x_1 x_2 - x_4^2)^{\tfrac{n}{d + 2 n}} \left(8 d (x_4^2-x_1 x_2)^2 - n (x_2^2 + x_1^2 (4 x_2^2 - 8 x_2 +1) - 2 x_4^2 + 8 x_2 x_4^2- 4 x_4^4 -8 x_1 (x_2^2 - x_4^2))\right)}{16 (x_4^2-x_1 x_2)^2}.$$
\end{corollary}

The Einstein metrics are precisely the critical points of the normalized scalar curvature. As mentioned in Section \ref{metric}, there may be either one or two non-diagonal Einstein metrics of the form  $g = (x_1, x_2, x_3, x_4)$ on $M = H \times  H/\Delta K$.

\begin{theorem}\label{stability}
If $H/K$ is any irreducible symmetric space, then the $H\times H$-invariant Einstein metric on the homogeneous space $M=H\times H/\Delta K$, $g_5=\left(\unm,\frac{3}{2},1,\unm\right)$ is non-degenerate, $G$-unstable of co-index $=1$ and if we assume $\kil_\kg=a\kil_\hg|_\kg$ (see \cite[Table 3]{HHK}) with $a>\tfrac{1}{2}$ the Einstein metric $g_3=(1,1,1,y_+)$ is non-degenerate, $G$-unstable with co-index $=2$.
\end{theorem}

\begin{remark}
  In particular, we obtain that $T_g\mca^{4}_1$ consists of $TT$-tensors (see \cite[Section 3.4]{stab-tres}).
\end{remark}

\begin{proof}
Consider the normalized scalar curvature $\scalar_N(g)$ given in Corollary \ref{snor}, $\left(\unm,\frac{3}{2},\unm\right)$ is a critical point of it, and the Hessian of this function at this critical point is:

\small{
$$H(g_5):=\begin{bmatrix} -\tfrac{2^{\tfrac{-n}{d + 2 n}} n (7 d + 23 n)}{2 (d + 2 n)}&-\tfrac{2^{\tfrac{-n}{d + 2 n}} n (-d + n)}{2 (d + 2 n)}&\tfrac{2^{\tfrac{-n}{d + 2 n}} n (d + 5 n)}{d + 2 n}\\
 -\tfrac{2^{\tfrac{-n}{d + 2 n}} n (-d + n)}{2 (d + 2 n)}&\tfrac{2^{\tfrac{-n}{d + 2 n}} n (d + n)}{2 (d + 2 n)}&-\tfrac{2^{\tfrac{-n}{d + 2 n}} n (d + n)}{d + 2 n} \\
\tfrac{2^{\tfrac{-n}{d + 2 n}} n (d + 5 n)}{d + 2 n} &\tfrac{2^{\tfrac{-n}{d + 2 n}} n (d + n)}{2 (d + 2 n)}& -\tfrac{2 (2^{\tfrac{-n}{d + 2 n}}) n^2}{d + 2 n}
\end{bmatrix}.
$$}

Note that the first and last elements of the diagonal are negative, and the second element is positive. This means that we can consider directions $e_1:=(1,0,0)$ and $e_2:=(0,1,0)$ such that $\langle H(g_5)e_1,e_1 \rangle <0$ and $\langle H(g_5)e_2,e_2 \rangle >0$, therefore the minimum eigenvalue of $H(g_5)$ is negative and the maximum is positive, thus $g_5$ is unstable. The determinant of this matrix is $2 (2^{-n})^{\tfrac{3}{d + 2 n}} n^3$, since it is always positive the co-index of $g_5$ is $1$ on $M$.

On the other hand, when $\kil_\kg=a\kil_\hg|_\kg$ and $a>\tfrac{1}{2}$ the equality $n=2d(1-a)$ holds and $g_3=(1,1,1,y_+)$, where $y_+=\frac{1}{2}\sqrt{\frac{2a-1}{2-a}}$ is a critical point of $\scalar_N$. The Hessian of this function at $(1,1,y_+)$ is:

\small{
$$H(g_3):=\begin{bmatrix}
\tfrac{4 \ C (a - 2)^2 (a - 1) d}{3 (3 - 2 a)^2 (4 a -5)}&
\tfrac{4 \ C (a - 2)(16a^3-53a^2+56a-19) d}{9 (3 - 2 a)^2 (4 a -5)}&
\tfrac{4 \ C \sqrt{\tfrac{1 - 2 a}{a-2}} (a - 2)^3(a - 1) d}{9 (3 - 2 a)^2 (4 a -5)}\\
 \tfrac{4 \ C (a - 2)(16a^3-53a^2+56a-19) d}{9 (3 - 2 a)^2 (4 a -5)}&\tfrac{4 \ C (a - 2)^2 (a - 1) d}{3 (3 - 2 a)^2 (4 a -5)}&\tfrac{4 \ C \sqrt{\tfrac{1 - 2 a}{a-2}}  (a - 2)^3(a - 1) d}{9 (3 - 2 a)^2 (4 a -5)} \\
\tfrac{4 \ C \sqrt{\tfrac{1 - 2 a}{a-2}}(a - 2)^3(a - 1) d}{9 (3 - 2 a)^2 (4 a -5)} &\tfrac{4 \ C \sqrt{\tfrac{1 - 2 a}{a-2}}  (a - 2)^3(a - 1) d}{9 (3 - 2 a)^2 (4 a -5)}&
\tfrac{16 \ C (a - 2)^2 (2 a^2 - 3 a + 1) d}{9 (3 - 2 a)^2 (4 a - 5)}
\end{bmatrix},
$$}

where $C:=\left(\tfrac{9 - 6 a}{8 - 4 a}\right)^{\scriptscriptstyle{\tfrac{-2 (a - 1) d}{5 d - 4 a d}}}>0$.

Note that the first two elements of the diagonal are equals and negative for all $0<a<1$ , and the third element is positive when $a>\tfrac{1}{2}$. As before, the minimum eigenvalue of $H(g_3)$ is negative and the maximum is positive, thus $g_3$ is unstable. Additionally, since the determinant of this matrix is given by
$$-\frac{256 C^3 (a - 2)^4 (a - 1)^3 (2 a -1) d^3}{
 243 (2 a-3)^5},$$
which is always negative for $a>\tfrac{1}{2}$, we can conclude that the co-index of $g_3$ is always $2$ on $M$.
\end{proof}

\begin{figure}[ht]
\begin{subfigure}{.5\textwidth}
    \centering
    \includegraphics[width=.95\linewidth]{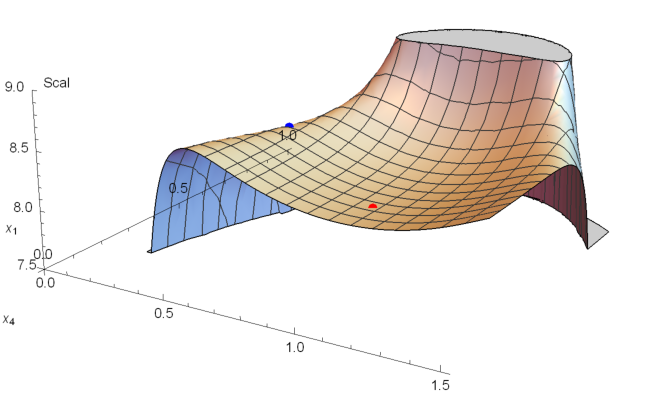}
\end{subfigure}%
\begin{subfigure}{.5\textwidth}
    \centering
    \includegraphics[width=.95\linewidth]{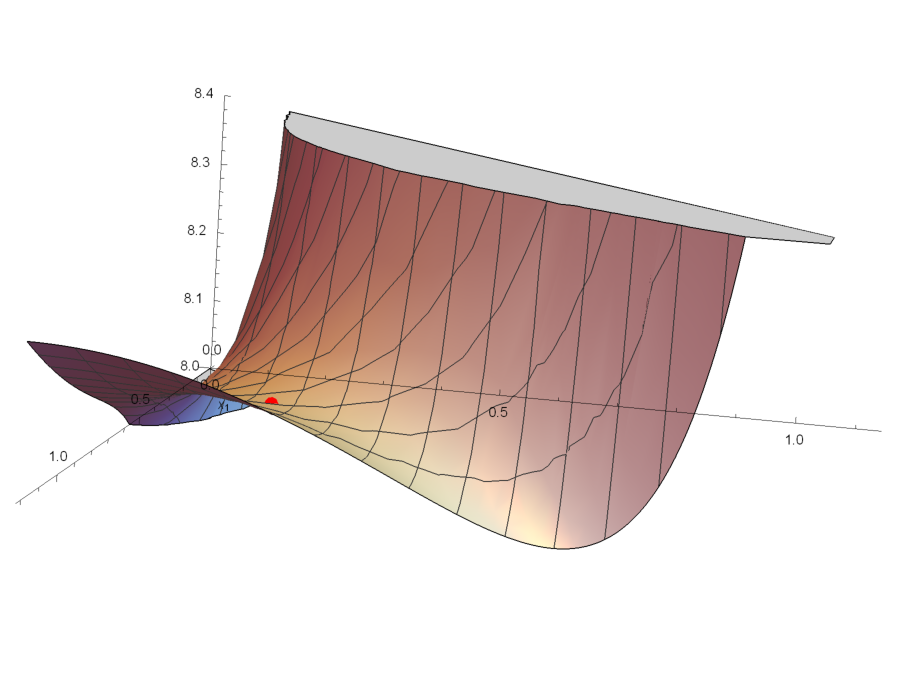}
\end{subfigure}
\caption{Graph of $\scalar_{N}$ restricted to the planes $x_1+x_2=2$ and $x_1=x_2$ on $M^{20}=\SU(4)\times \SU(4) /\Delta \Spe(2)$, which admits two non-diagonal Einstein metrics up to isometry, {\color{blue} $g_5$} (in blue) and {\color{red} $g_3$} (in red).}
\label{figu}
\end{figure}

If we consider $M^{20}=\SU(4)\times \SU(4) /\Delta \Spe(2)$, it was proved that there are two distinct non-diagonal Einstein metrics on $M$ of the form $g=(x_1,x_2,x_3,x_4)$ ($a=\tfrac{3}{4}$, see \cite[Table 3]{HHK}). On Figure \ref{figu}, they were illustrated as critical points of the normalized scalar curvature given by Corollary \ref{snor}. While the metric $g_3$ appears to be a minimum when projected onto the plane $x_1+x_2=2$, it is actually found to have a co-index of $2$ when considered in the projection onto the plane where $x_1=x_2$.

\end{document}